\documentclass[12pt]{article}
\usepackage{amsmath}
\usepackage{amssymb}
\usepackage{vector}
\def\eps{\varepsilon}
\font\tencmmib=cmmib10 \skewchar\tencmmib '60
\newfam\cmmibfam
\textfont\cmmibfam=\tencmmib

\def\bbox{\quad\hbox{\vrule \vbox{\hrule \vskip2pt \hbox{\hskip2pt
\vbox{\hsize=1pt}\hskip2pt} \vskip2pt\hrule}\vrule}}
\def\lessim{\ \lower4pt\hbox{$
\buildrel{\displaystyle <}\over\sim$}\ }
\def\gessim{\ \lower4pt\hbox{$\buildrel{\displaystyle >}
\over\sim$}\ }

\def\eps{{\varepsilon}}

\def\la{{\Bigl\langle}}
\def\ra{{\Bigr\rangle}}

\def\qed{\hfill\break\rightline{$\bbox$}}
\newcommand{\e}{\mathbb{E}}
\newcommand{\p}{\mathbb{P}}
\newcommand{\Reals}{\mathbb{R}}

\newtheorem{lemma}{Lemma}
\newtheorem{theorem}{Theorem}
\newtheorem{corollary}{Corollary}
\makeatletter
\@addtoreset{equation}{section}

\makeatother

\font\tencmmib=cmmib10 \skewchar\tencmmib '60
\newfam\cmmibfam
\textfont\cmmibfam=\tencmmib

\def\bbox{\quad\hbox{\vrule \vbox{\hrule \vskip2pt \hbox{\hskip2pt
\vbox{\hsize=1pt}\hskip2pt} \vskip2pt\hrule}\vrule}}
\def\lessim{\ \lower4pt\hbox{$
\buildrel{\displaystyle <}\over\sim$}\ }
\def\gessim{\ \lower4pt\hbox{$\buildrel{\displaystyle >}
\over\sim$}\ }

\def\AA{{\cal A} }

\def\eps{\varepsilon}

\def\go0{\to 0}

\def\la{\langle}

\def\leftitem#1{\item{\hbox to\parindent{\enspace#1\hfill}}}

\def\qed{\hfill\break\rightline{$\bbox$}}

\def\ra{\rangle}

\def\sg{\sigma}

\def\sg2{\sigma^2}

\def\__{_{\infty}}

\begin{document}

\title{
A note on Talagrand's positivity principle.}

\author{ 
{\it Dmitry Panchenko}\thanks{
Department of Mathematics: 
Texas A\&M University, College Station, TX and
Massachusetts Institute of Technology, Cambridge, MA.
Email: panchenk@math.tamu.edu. Partially supported by NSF grant.
}
}
\date{}

\maketitle

\begin{abstract}
Talagrand's positivity principle, Section 6.6 in \cite{SG}, states that one
can slightly perturb a Hamiltonian in the Sherrington-Kirkpatrick model
in such a way that the overlap of two configurations under 
the perturbed Gibbs' measure will become typically
nonnegative.  In this note we observe that abstracting from the setting
of the SK model only improves the result and does not require any modifications
in Talagrand's argument.  In this version, for example, positivity principle 
immediately applies to the setting in \cite{T2}. Also, abstracting from the SK
model improves the conditions in the Ghirlanda-Guerra identities and as a consequence
results in a perturbation of smaller order necessary to ensure positivity of
the overlap.
\end{abstract}
\vspace{0.5cm}



\section{Introduction.}

Let us consider a unit sphere $S=\{z\in \Reals^N: |z|=1\}$ on euclidean space $\Reals^N$
and let $\nu$ be a probability measure on $S.$ Given a measurable function $g: S\to\Reals$
let us define a probability measure $\nu_g$ on $S$ by a change of density 
$$
\frac{d \nu_g}{d \nu} = \frac{e^{g(z)}}{\int e^{g(z)} d\nu(z)}.
$$
We assume that the denominator on the right hand side is integrable.
Let us make a specific choice of $g(z)$ given by
\begin{equation}
g(z) = v \sum_{p\geq 1} 2^{-p} x_p\, g_p(z)
\,\,\,\mbox{ for }\,\,\,
g_p(z)= \sum_{1\leq i_1,\ldots, i_p\leq N}g_{i_1,\ldots,i_p} z_{i_1}\cdots z_{i_p},
\label{g}
\end{equation}
where $v\geq 0,$ $x_p$ are i.i.d. random variables uniform on $[0,1]$ and
$g_{i_1,\ldots,i_p}$ are i.i.d. standard Gaussian for all $i_1,\ldots, i_p$
and all $p\geq 1.$ 

Given a function $f$ on $S^n$ let us denote by $\la f\ra$ its average with respect
to measure $\nu_g^{\otimes n}.$
Let us denote by $\e_g$ the expectation with respect to Gaussian random variables
and by $\e_x$ the expectation with respect to uniform random variables in the 
definition of $g(z)$ in (\ref{g}).
The following is the main result of this note.

\begin{theorem}\label{Th1}{\rm (Talagrand's positivity principle)}
For any $\eps>0$ there exists large enough $v\geq 1$ in (\ref{g}) such that 
\begin{equation}
\e \nu_g^{\otimes 2}\bigl\{z^1\cdot z^2 \leq -\eps\bigr\} \leq \eps.
\label{pos}
\end{equation}
The choice of $v$ does not depend on $N$ and $\nu.$
\end{theorem}

This means that one can define a random perturbation $\nu_g$ of an arbitrary
measure $\nu$ such that the scalar product $z^1\cdot z^2$ of two vectors drawn 
independently from distribution $\nu_g$ will be typically nonnegative. 
This result was proved in Section 6.6 \cite{SG} in the setting of 
the Sherrington-Kirkpatrick model where $\nu$ was a random Gibbs' measure 
and in which case the expectation on the left hand side of (\ref{pos}) 
was also in the randomness of $\nu$. The main ingredient of the proof was
the Ghirlanda-Guerra identities that are typically given in the setting of the SK model as well.

The main contribution of this note is an observation that
abstracting ourselves from the setting of the SK model results in some qualitative
improvements of the positivity principle of Theorem \ref{Th1}. 
First of all, we notice that Talagrand's proof in \cite{SG} requires no modifications 
in order to prove a positivity principle that holds uniformly over $\nu$ rather than on
average over a random Gibbs' measure in the SK model.
This observation, for example, implies the result in \cite{T2} without any additional work 
as will be shown in Example below. 
Another important qualitative improvement is the fact that the choice of $v$ in (\ref{pos}) is 
independent of $N.$ In \cite{SG}, one needed $v\gg N^{1/4}$ - a condition that appears in the proof of 
Ghirlanda-Guerra identities due to the fact that one controls random Gibbs' measure from the very beginning. 
We will show below that the supremum of $g(z)$ on the sphere is of order $v\sqrt{N}$ which means
that one can perturb any measure on $S$ by a change of density of order $\exp v\sqrt{N}$ and
force the scalar product $z^1\cdot z^2$ to be essentially nonnegative.

{\it Example}
{\it (Positivity in Guerra's replica symmetry breaking bound, \cite{T2}).} 
The main result in \cite{T2} states that Guerra's replica symmetry breaking bound \cite{Guerra}
applies to odd $p$-spin interactions as well. The proof utilizes the Aizenman-Sims-Starr
version of Guerra's interpolation \cite{ASS} and a positivity principle that requires 
a concentration inequality for the free energy along the interpolation. 
We observe that this positivity principle follows directly from Theorem \ref{Th1}.
Let $\AA$ be a countable set and let $(w_{\alpha})$ be a probability function on $\AA$ 
such that 
$
w_{\alpha}\geq 0,\,\,\,\sum_{\alpha\in\AA}w_{\alpha}=1.
$ 
Let $H(z,\alpha)$ be a function on $\Omega=\Sigma\times\AA$ 
for some finite subset $\Sigma$ of  $S.$  
Let us consider a probability measure on $\Omega$ given by
$$
\mu\{(z,\alpha)\}\sim w_{\alpha}\exp\bigl( H(z,\alpha) + g(z)\bigr)
$$
where $g(z)$ is defined in (\ref{g}). Then its marginal on $\Sigma$ is
equal to $\nu_g$ if we define $\nu$ by
$
\nu\{z\} \sim \sum w_{\alpha}\exp H(z,\alpha).
$
Therefore,
$$
\mu^{\otimes 2}\bigl\{z^1\cdot z^2\leq -\eps\bigr\}
=
\nu_g^{\otimes 2}\bigl\{z^1\cdot z^2\leq -\eps\bigr\}.
$$
By Theorem \ref{Th1}, for large enough $v>0,$
$$
\e \mu^{\otimes 2}\bigl\{z^1\cdot z^2\leq -\eps\bigr\}\leq \eps
$$
and this inequality holds uniformly over all choices of $H$ and $w.$
Therefore, we can average over arbitrary random distribution of $H$ and $w.$
In particular, in \cite{T2}, $(w_{\alpha})$ was a random Derrida-Ruelle process on
$\AA=\mathbb{N}^k$ and $H$ was Guerra's interpolating Hamiltonian  
in the form of Aizenman-Sims-Starr.
\qed

{\it General remarks.}
(i) As we will see below, Talagrand's proof of positivity principle uses very
deep information about a joint distribution of $z^1\cdot z^2$ for $2\leq l\leq n$
under measure $\e_g \nu_g^{\otimes n}$ for a typical realization of $(x_p).$
However, it is not clear whether this deep information is really necessary to prove
positivity and if some simpler argument would not suffice. For example, if we consider
only a first order Gaussian term in (\ref{g}) and define $g'(z)=v\sum_{i\leq N} g_i z_i$  
then measure $\nu_{g'}$ under a change of density proportional to $e^{g'(z)}$ would
favor a random direction $g=(g_1,\ldots,g_N)$ and it is conceivable that for large
enough $v$ independent of $\nu$ and $N$ two independent vectors $z^1$ and $z^2$ from
this measure would typically "point in the same direction", i.e.
$$
\e \nu_{g'}^{\otimes 2}\{z^1\cdot z^2 \leq -\eps\}
=
\e\, \frac{\int I(z^1\cdot z^2\leq -\eps) e^{v(g\cdot z^1 + g\cdot z^2)} d\nu(z^1)d\nu(z^2)}
{\bigl(\int e^{v g\cdot z} d\nu(z)\bigr)^2}
\leq \eps.
$$ 
In fact, even a weaker result with $v=o(\sqrt{N})$ 
would be sufficient for applications as in \cite{T2}. 

(ii) Theorem \ref{Th1} implies the following non-random version of positivity
principle.

\begin{corollary}\label{Cor1}
For any $\eps>0$ there exists $v>0$ large enough such that the following holds.
For any distribution $Q$ on the set of measures on the sphere $S$ 
there exists a (non-random) function $g(z)$ such that for some absolute constant $L,$
$$
\sup_{z\in S}|g(z)| \leq L v\sqrt{N}
\,\,\mbox{ and }\,\,
\int \nu_g^{\otimes 2}\{z^1\cdot z^2 \leq -\eps\} dQ(\nu)\leq \eps.
$$
\end{corollary}
It would be of interest to prove this result directly and not as a corollary
of Theorem \ref{Th1}. Then, by Hahn-Banach theorem, one can find a distribution
$P$ on the set of functions $\{\sup_{z\in S}|g(z)| \leq Lv\sqrt{N}\}$ such that
for all probability measures $\nu$ on $S$,
$$
\int \nu_g^{\otimes 2}\{z^1\cdot z^2 \leq -\eps\} dP(g) \leq \eps.
$$
This would give another proof of Theorem \ref{Th1} with non constructive
description of $P.$
\qed

{\it Sketch of the proof of positivity principle.}
The main ingredient in Talagrand's proof is the 
extended Ghirlanda-Guerra identities (Section 6.4 in \cite{SG})
that state that if we sample 
$(z^1,\ldots,z^{n+1})$ from a measure $\e_g\nu_g^{\otimes (n+1)}$ 
(which is a mixture of product measures of $\nu_g$ over the randomness
of Gaussian r.v.s)
then for a typical realization of $(x_p)$ the scalar product $z^1\cdot z^{n+1}$ with probability
$1/n$ is independent of $z^1\cdot z^l, 2\leq l\leq n,$ and with probabilities
$1/n$ it is equal to one of them. More precisely, the following holds.

\begin{theorem}\label{Th2} {\rm (Ghirlanda-Guerra identities) }
For any measurable function $f$ on $S^n$ such that $|f|\leq 1$ and
for any continuous function $\psi$ on $[-1,1]$ we have
\begin{equation}
\e_x
\Bigl|
\e_g\la f \psi(z^1\cdot z^{n+1})\ra 
-\frac{1}{n}\e_g\la f\ra \e_g\la \psi(z^1\cdot z^2)\ra
-\frac{1}{n}\sum_{2\leq l\leq n} 
\e_g\la f \psi(z^1\cdot z^l)\ra 
\Bigr|\leq \delta(\psi,n,v)
\label{GG}
\end{equation}
where $\delta(\psi,n,v)\to 0$ as $v\to\infty$ and $\delta$ does not
depend on $N$ and $\nu.$
\end{theorem}

The main idea of Talagrand's proof can be shortly described as follows.
Suppose that we sample $z^1,\ldots,z^n$ independently from any measure on $S.$
Then the event that all $z^1\cdot z^l\leq -\eps$ simultaneously is very unlikely
and its probability is or order $1/(n\eps).$ The bound on this probability is uniform over
all measures and therefore can be averaged over some distribution on measures
and holds, for example, for $\e_g \nu_g^{\otimes n}.$ On the other hand,
by Ghirlanda-Guerra identities under measure $\e_g \nu_g^{\otimes n}$
 the events $\{z^1\cdot z^l \leq -\eps\}$ are strongly
correlated due to the fact that with some prescribed probabilities $z^1\cdot z^l$
can be equal. As a consequence, the simultaneous occurrence of these events will
have probability of order $a/n^{1-a}$ where $a$ is the probability of one of them.
But since this quantity is bounded by $1/(n\eps),$ taking $n$ large enough shows that
$a$ should be small enough. 
\qed

\section{Proofs.}

The proof of Theorem \ref{Th2} follows exactly the same argument as
the proof of Ghirlanda-Guerra identities in \cite{GG} or in \cite{SG}.
Since we consider a fixed measure $\nu,$ we do not need to control a
fluctuations of a random Hamiltonian and, as a result, 
we get a better condition on $v.$
The main part of the proof is the following lemma.

\begin{lemma}\label{L1}
For any $p\geq 1$ there exists a constant $L_p$ that depends on $p$ only such that
\begin{equation}
\e \bigl\la \bigl|g_p(z) -\e_g\bigl\la g_p(z)\bigr\ra \bigr|\bigr\ra \leq L_p\sqrt{v}.
\label{GGbasic}
\end{equation}
\end{lemma}
{\it Proof.}
Let us fix $(x_p)_{p\geq 1}$ in the definition of $g(z)$ and until the end of the proof
let $\e$ denote the expectation in Gaussian random variables only. Define
$$
\theta = \log \int e^{g(z)} d\nu(z),\,\,\,
\psi = \e \theta.
$$ 
Given $p\geq 1,$ let us think of $\theta$ and $\psi$ as functions of $x=x_p$ only
and define $v_p = v 2^{-p}.$ 
Then $\theta'(x)=v_p \bigl\la g_p(z)\bigr\ra$ and
\begin{equation}
\psi'(x)=v_p \e\bigl\la g_p(z)\bigr\ra
=  v_p^2 x\bigl(1-\e\la (z^1\cdot z^2)^p\ra\bigr)
\leq v_p^2 x.
\label{der}
\end{equation}
Since
$$
\psi''(x)= v_p^2 \e \bigl( \bigl\la g_p(z)^2 \bigr\ra - \bigl\la g_p(z)\bigr\ra^2\bigr)
= v_p^2 \e\bigl\la \bigl(g_p(z) - \la g_p(z)\ra\bigr)^2\bigr\ra.
$$
we have
$$
v_p^2 
\int_{0}^{1} \e\bigl\la \bigl(g_p(z) - \la g_p(z)\ra\bigr)^2\bigr\ra dx
=\psi'(1)-\psi'(0)\leq v_p^2
$$
and by Cauchy inequality
\begin{equation}
\int_{0}^{1} \e\bigl\la \bigl|g_p(z) - \la g_p(z)\ra\bigr| \bigr\ra dx \leq 1.
\label{step2}
\end{equation}
To prove (\ref{GGbasic}) is remains to approximate $\la g_p(z)\ra$ by $\e\la g_p(z)\ra$ 
and to achieve that we will use a simple consequence of convexity of $\theta$ and $\psi$
given in the inequality (\ref{conv}) below. Since 
$$
\e g(z)^2 = v^2 \sum_{p\geq 1} 2^{-2p} x_p^2 \leq 2 v^2
$$
if all $|x_p| \leq 2,$ we can apply a well-known Gaussian concentration inequality, 
for completeness given in Lemma \ref{Lextra4} in Appendix A, to get
$$
\e |\theta(x) -\psi(x)| \leq 4 v 
\,\,\,\mbox{ and }\,\,\,
\e |\theta(x\pm y) -\psi(x\pm y)| \leq 4 v.
$$ 
Below we will choose $|y|\leq 1$ so that $|x\pm y|\leq 2.$
({\it Remark.} At the same step in the setting of the SK model one also needs to control
a random Hamiltonian which produces another term of order $\sqrt{N}$
and this results in unnecessary condition on $v$ in the positivity principle.)
Inequality (\ref{conv}) then implies that
$$
\e|\theta'(x)-\psi'(x)| 
\leq
\psi'(x+y) -\psi'(x-y) 
+ \frac{12 v}{y}.
$$
Since
$$
\int_{0}^{1} (\psi'(x+y) - \psi'(x-y))dx 
= 
\psi(1+y)-\psi(1-y)-\psi(y)+\psi(-y)
\leq 
6 v_p^2 y
$$ 
using (\ref{der}) and $|y|\leq 1$ we get
$$
\int_{0}^{1} \e|\theta'(x)-\psi'(x)| dx
\leq
6 v_p^2 y + \frac{12 v}{y}\leq L_p v^{3/2}
$$
if we take $y=v^{-1/2}\leq 1.$ Using explicit expressions for
$\theta'$ and $\psi'$ we finally get
\begin{equation}
\int_{0}^{1} \e|\la g_p(z)\ra-\e\la g_p(z)\ra| dx
\leq L_p \sqrt{v}.
\label{step1}
\end{equation}
Together with (\ref{step2}) this gives
$$
\int_{0}^{1} \e \bigl\la | g_p(z)-\e\la g_p(z)\ra |\bigr\ra dx
\leq L_p \sqrt{v}.
$$
We now recall that $\e$ was expectation $\e_g$ with respect to Gaussian 
random variables only. Also, integral over $x=x_p\in [0,1]$ is nothing but
expectation with respect to $x_p.$ Therefore, averaging over all remaining
$x_p$ finishes the proof.
\qed

The following inequality was used in the previous proof and it 
quantifies the fact that if two convex functions are close
to each other then their derivatives are also close.

\begin{lemma}
If $\theta(x)$ and $\psi(x)$ are convex differentiable functions then
\begin{eqnarray}
|\theta'(x)-\psi'(x)|
&\leq& 
\psi'(x+y) -\psi'(x-y) 
\label{conv}
\\
&+& 
\frac{1}{y}\bigl(|\psi(x+y)-\theta(x+y)|+|\psi(x-y)-\theta(x-y)| 
+ |\psi(x) - \theta(x)|\bigr).
\nonumber
\end{eqnarray}
\end{lemma}
{\it Proof.}
By convexity, for any $y>0$
$$
\frac{\theta(x)-\theta(x-y)}{y}\leq \theta'(x) 
\leq 
\frac{\theta(x+y)-\theta(x)}{y}
$$
and
$$
\frac{\psi(x)-\psi(x-y)}{y}\leq \psi'(x) 
\leq 
\frac{\psi(x+y)-\psi(x)}{y}.
$$
If we define
$$
U = 
\frac{1}{y}\Bigl(|\psi(x+y)-\theta(x+y)|+|\psi(x-y)-\theta(x-y)| 
+ |\psi(x) - \theta(x)|\Bigr)
$$
then the above inequalities imply
\begin{eqnarray*}
\theta'(x)
&\leq&
\frac{\theta(x+y)-\theta(x)}{y}
\leq
\frac{\psi(x+y)-\psi(x)}{y} +U
\\
&=&
\frac{\psi(x)-\psi(x-y)}{y}
+
\frac{\psi(x+y)+\psi(x-y)-2\psi(x)}{y} +U
\\
&\leq&
\psi'(x) + \psi'(x+y) -\psi'(x-y) +U.
\end{eqnarray*}
Similarly,
$$
\theta'(x)\geq
\psi'(x) - (\psi'(x+y) -\psi'(x-y))-U.
$$
Combining two inequalities finishes the proof.
\qed

{\it Proof of Theorem \ref{Th2}.}
Since $|f|\leq 1$ we can write
$$
\Bigl|\e_g\bigl\la f g_p(z) \bigr\ra 
- \e_g\la f \ra \e\la g_p(z)\ra \Bigr|
\leq
\e_g\bigl\la \bigl| g_p(z) - \e_g \la g_p(z)\ra\bigr|\bigr\ra. 
$$
By Gaussian integration by parts, the left hand side is equal to 
$n v 2^{-p}x_p \phi_p$ where
$$
\phi_p=\Bigl|
\e_g\bigl\la f (z^1\cdot z^{n+1})^p\bigr\ra 
-\frac{1}{n}\e_g\bigl\la f\ra\e_g\la (z^1\cdot z^2)^p\bigr\ra
-\frac{1}{n}\sum_{2\leq l\leq n} 
\e_g\bigl\la f (z^1\cdot z^l)^p\bigr\ra 
\Bigr|.
$$
Lemma \ref{L1} then implies $nv2^{-p} \e_x x_p \phi_p \leq L_p \sqrt{v}$
and, thus, 
$$
\e_x x_p \phi_p \leq Lv^{-1/2}
$$ 
for some constant $L$ that depends only on $n$ and $p$.
Since $\phi_p\leq 2,$ for any $x_0\in(0,1),$
$$
\e_x \phi_p \leq 2x_0 + \frac{1}{x_0} \e_x x_p\phi_p
\leq 2x_0 + \frac{Lv^{-1/2}}{x_0}
$$
and minimizing over $x_0$ we get $\e_x \phi_p(x_p)  \leq Lv^{-1/4}.$
Since any continuous function $\psi$ on $[-1,1]$ can be approximated
by polynomials this, obviously, implies the result.
\qed

{\it Proof of Theorem \ref{Th1}.}
{\it Step 1.} First we use Ghirlanda-Guerra identities to give a lower
bound on probability that all $z^1\cdot z^l\leq -\eps$ for $2\leq l\leq n.$
In order to use Theorem \ref{Th2} it will be convenient to denote
by $\delta$ any quantity that depends on $(x_p)$ and such that
$\e_x |\delta|$ does not depend on $\nu$ and $N$ and $\e_x |\delta|\to 0$
as $v\to \infty.$ Then (\ref{GG}) can be written as
$$
\e_g\la f \psi(z^1\cdot z^{n+1})\ra 
-\frac{1}{n}\e_g\la f\ra \e_g\la \psi(z^1\cdot z^2)\ra
-\frac{1}{n}\sum_{2\leq l\leq n} 
\e_g\la f \psi(z^1\cdot z^l)\ra 
=\delta.
$$
Even though a function $\psi$ is assumed to be continuous, let us use
this result formally for $\psi(x)=I(x\leq -\eps).$ 
The argument can be easily modified by using continuous approximations of 
the indicator function $\psi.$ Let
$$
f_n = \prod_{2\leq l\leq n} \psi(z^1\cdot z^l) 
= I(z^1\cdot z^l\leq -\eps \mbox{ for } 2\leq l\leq n)
$$
and let 
$$
a= \e_g\la \psi(z^1\cdot z^2)\ra = \e_g \nu_g^{\otimes 2}\{z^1\cdot z^2\leq-\eps\}.
$$
Then by Ghirlanda-Guerra identities
$$
\e_g\la f_{n+1} \ra = 
\e_g\la f_n \psi(z^1\cdot z^{n+1})\ra  
=
\frac{1}{n}\e_g\la f_n\ra \e_g\la \psi(z^1\cdot z^2)\ra
+\frac{1}{n}\sum_{2\leq l\leq n} 
\e_g\la f_n \psi(z^1\cdot z^l)\ra 
+\delta
$$
and since $f_n \psi(z^1\cdot z^l)=f_n$ we get
$$
\e_g\la f_{n+1} \ra = \frac{n-1+a}{n}\, \e_g\la f_{n} \ra +\delta.
$$
By induction,
\begin{equation}
\e_g\la f_n \ra = a \prod_{2\leq l\leq n-1}\frac{l-1+a}{l} + \delta
\geq
\frac{a}{Ln^{1-a}}+\delta
\label{GGpos}
\end{equation}
where the last inequality follows from a simple estimate for $l\geq 2$
$$
\frac{l-1+a}{l}=1-\frac{1-a}{l}\geq
\exp\Bigl(
-\frac{1-a}{l}-\frac{1}{l^2}
\Bigr).
$$
{\it Step 2.} On the other hand, we will show that $\la f_n \ra$ is of order $1/(n\eps)$
and to emphasize the fact that this is true for any measure we now simply write $G$
instead of $\nu_g.$ 
If $z^1,\ldots,z^n$ are i.i.d. from distribution $G$ then
\begin{equation}
\la f_n \ra 
=
\Bigl\la G\bigl\{z^2 : z^1\cdot z^2 \leq -\eps\bigr\}^{n-1}\Bigr\ra 
\leq 
G(U) + \gamma^{n-1}
\label{GU}
\end{equation}
where, given $0<\gamma<1,$ we defined a set
$$
U=\Bigl\{
z^1\in S : G\{z^2 : z^1\cdot z^2\leq -\eps\} \geq \gamma
\Bigr\}.
$$
We would like to show that if $\gamma$ is close to $1$ then $G(U)$ is small.
This follows from the fact that the average of $z^1\cdot z^2$ is nonnegative
with respect to any product measure,
$$
\la z^1\cdot z^2\ra 
=
\sum_{i\leq N} \bigl\la z_i^1 z_i^2\bigr\ra
=
\sum_{i\leq N} \la z_i\ra^2 \geq 0.
$$
Since $z^1\cdot z^2 \leq 1,$
$$
0\leq \la z^1\cdot z^2\ra 
\leq 
-\eps G^{\otimes 2}\{z^1\cdot z^2 \leq -\eps\} 
+ G^{\otimes 2}\{z^1\cdot z^2> -\eps\}
=
-\eps +(1+\eps)G^{\otimes 2}\{z^1\cdot z^2> -\eps\}
$$
and
\begin{equation}
\eps\leq (1+\eps)G^{\otimes 2}\{z^1\cdot z^2> -\eps\}.
\label{pos1}
\end{equation}
Let us define a conditional distribution 
$G_U(C)=G(UC)/G(U).$
If $z^1\in U$ then by definition 
$G\{z^2 : z^1\cdot z^2> -\eps\} \leq 1-\gamma$
and
$$
G_U\{z^2: z^1\cdot z^2> -\eps\}
\leq
\frac{1-\gamma}{G(U)}.
$$
Since $G_U$ is concentrated on $U,$ 
$$
G_U^{\otimes 2}\{z^1\cdot z^2> -\eps\}\leq \frac{1-\gamma}{G(U)}.
$$
Using (\ref{pos1}) for $G_U$ instead of $G$ we get
$$
\eps \leq (1+\eps)G_U^{\otimes 2}\{z^1\cdot z^2> -\eps\} \leq \frac{2(1-\gamma)}{G(U)}
$$
and, therefore, $G(U)\leq 2(1-\gamma)/\eps.$ Then, (\ref{GU}) implies that
\begin{eqnarray*}
&&
\la f_n \ra 
\leq 
\frac{2(1-\gamma)}{\eps} + \gamma^{n-1}
\leq
\frac{2(1-\gamma)}{\eps} + e^{-(n-1)(1-\gamma)}
\end{eqnarray*}
and we can minimize this bound over $0<\gamma<1$ to get
$$
\la f_n \ra 
\leq 
\frac{L}{n\eps}\log\frac{n\eps}{L}.
$$
{\it Step 3.} Together with (\ref{GGpos}) this implies
$$
\frac{a}{Ln^{1-a}}\leq
\frac{L}{n\eps}\log\frac{n\eps}{L}
+ \delta.
$$
By definition of $\delta$ this implies that
$$
\e_x a n^{a } \leq \frac{L}{\eps}\log n\eps + \e_x|\delta|
$$
where $\e_x|\delta|\to 0$ as $v\to\infty.$
For any $a_0>0$ we can write
\begin{eqnarray*}
&&
\e_x a 
\leq
a_0 + \e_x a I(a\geq a_0) 
\leq
a_0 + \frac{1}{n^{a_0}} \e_x a n^{a}
\leq
a_0+
\frac{1}{n^{a_0}}
\Bigl(
\frac{L}{\eps}\log n\eps + \e_x|\delta|
\Bigr).
\end{eqnarray*}
The right hand side goes to zero if we let $v\to\infty, n\to\infty, a_0\to 0$
which means that $\e_x a \to 0$ as $v\to\infty$ and this finishes the proof.
\qed

{\it Proof of Corollary \ref{Cor1}.}
Since (\ref{pos}) is uniform in $\nu,$ by Chebyshev's inequality
\begin{equation}
\p\Bigl(
\int \nu_g^{\otimes 2}\{z^1\cdot z^2 \leq -\eps\} dQ(\nu)
\geq 4\eps
\Bigr)\leq \frac{1}{4}.
\label{half}
\end{equation}
Next, we will show that for some absolute constant $L,$
\begin{equation}
\p\Bigl(
\sup_{z\in S}|g(z)|
\geq 
Lv\sqrt{N} \Bigr)
\leq \frac{1}{4}.
\label{halftwo}
\end{equation}
Indeed, conditionally on $(x_p)$ the process $g(z)$ is Gaussian with covariance
$$
\e_g g(z^1) g(z^2) = v^2 \xi(z^1\cdot z^2)
\,\,\mbox{ where }\,\,
\xi(s) = \sum_{p\geq 1} 2^{-2p} x_p^2 s^p.
$$
Since the first two derivatives of $\xi(s)$ are bounded uniformly over $(x_p)\in [0,1]^{\infty},$
\begin{eqnarray*}
\e_g (g(z^1)-g(z^2))^2 
&=& 
v^2 
\bigl(\xi(z^1\cdot z^1)+\xi(z^2\cdot z^2) - 2\xi(z^1\cdot z^2)
\bigr)
\\
&\leq& 
L v^2 |z^1-z^2|^2 = Lv^2 \e(\eta(z^1) - \eta(z^2))^2
\end{eqnarray*}
where $\eta$ is the canonical Gaussian process $\eta(z) = \sum g_i z_i.$
By a consequence of Slepian's inequality (see Corollary 3.14 in \cite{LT})
$$
\e_g \sup_{z\in S} g(z) \leq Lv\, \e \sup_{z\in S} \eta(z)
\leq L v\sqrt{N} 
$$
and, by symmetry, $\e_g \sup_{z\in S} |g(z)| \leq Lv\sqrt{N}$
which implies (\ref{halftwo}). (\ref{half}) and (\ref{halftwo}) imply that
the event
$$
\Bigl\{
\int \nu_g^{\otimes 2}\{z^1\cdot z^2 \leq -\eps\} dQ(\nu)
\leq 4\eps
\Bigr\}\bigcap
\Bigl\{
\sup_{z\in S}|g(z)|
\leq 
Lv\sqrt{N}\Bigr\}
$$
is not empty and this finishes the proof.
\qed

\appendix

\section{A gaussian concentration inequality.}

For completeness, we give a proof of the Gaussian concentration
inequality which can be found, for example, in \cite{SKgen}. 

\begin{lemma}\label{Lextra4}
Let $\nu$ be a finite measure and $g(z)$ be a Gaussian process on $\Reals^N$ such that 
$
\e g(z)^2 \leq a 
$
for $z$ in the support of measure $\nu.$ If
$$
X=\log\int \exp g(z) d\nu(z)
$$
then $\e(X-\e X)^2 \leq 8 a.$
\end{lemma}
\textbf{Proof.}
Let $g^1$ and $g^2$ be two independent copies of 
$g.$ For $t\in [0,1]$ and $j=1,2$ we define
$$
g_t^j(z) = \sqrt{t} g^j(z)+\sqrt{1-t}g(z)
\,\,\,\mbox{ and }\,\,\,
g_t(z^1,z^2)=g_t^1(z^1)+g_t^2(z^2).
$$
Let
\begin{equation}
F_j = \log \int \exp g_t^j(z) d\nu(z) 
\label{Pfi}
\end{equation}
and let $F=F_1+F_2.$ For $s\geq 0,$ let
$\varphi(t) = \e \exp s(F_2 -F_1).$
We can write
$$
\varphi'(t) 
=
s\e \exp s(F_2 -F_1)
\sum_{j\leq 2}(-1)^{j+1}\exp(-F_j)
\int \frac{\partial g_t^j(z)}{\partial t} 
\exp g_t^j(z) d\nu(z).
$$ 
If $\zeta(z^1,z^2):=\e g(z^1) g(z^2)$ then
$$
2\e \frac{\partial g_t^j(z^1)}{\partial t} g_t^j(z^2)
= \zeta(z^1,z^2) - \zeta(z^1,z^2)=0
$$
and
$$
2\e \frac{\partial g_t^1(z^1)}{\partial t} g_t^2(z^2)
=-\zeta(z^1,z^2).
$$
Therefore, Gaussian integration by parts gives
$$ 
\varphi'(t) =
-s^2 
\e \exp s(F_1 - F_2) \exp(-F)
\int \zeta(z^1,z^2)\exp g_t(z^1,z^2)
d\nu(z^1)d\nu(z^2)
$$
and since by assumption $|\zeta(z^1,z^2)|\leq a$ on the support of
measure $\nu^{\otimes 2}$ we get
$$
\varphi'(t) \leq a s^2 
\e \exp s(F_1 - F_2) \exp(-F)
\int \exp g_t(z^1,z^2) d\nu(z^1) d\nu(z^2)
= a s^2 \varphi(t).
$$
By construction $\varphi(0)=1$ and the above inequality implies that
$\varphi(1)\leq \exp as^2.$ 
On the other hand, by construction, $\varphi(1)=\e \exp s(X-X'),$
where $X'$ is an independent copy of $X.$ Thus,
$$
\e \exp s(X-X')\leq \exp as^2.
$$
By Jensen's inequality, 
$\e \exp s(X-\e X)\leq \exp as^2$ and by Markov's inequality 
$$
\p\bigl( X-\e X\geq t\bigr)\leq \inf_{s> 0} \exp(as^2 - st)
=\exp\Bigl(-\frac{t^2}{4a}\Bigr).
$$
Obviously, a similar inequality can be written for $\e X -X$ and, 
therefore,
$$
\p\bigl( |X-\e X|\geq t\bigr)\leq 2\exp\Bigl(-\frac{t^2}{4a}\Bigr).
$$
The result follows.
\qed


\begin{thebibliography}{99}

\bibitem{ASS} Aizenman, M., Sims, R., Starr, S. (2003)
An extended variational principle for the SK spin-glass model.
{\it Phys. Rev. B}, \textbf{68}, 214403. 

\bibitem{GG} Ghirlanda, S., Guerra, F. (1998) 
General properties of overlap probability distributions in disordered spin systems. 
Towards Parisi ultrametricity.  {\it J. Phys. A}  \textbf{31}, no. 46, 9149-9155.

\bibitem{Guerra} Guerra, F. (2003) Broken replica 
symmetry bounds in the mean field spin glass model. 
{\it Comm. Math. Phys.} {\bf 233}, no. 1, 1-12.

\bibitem{LT}
Ledoux, M., Talagrand, M. (1991)
Probability in Banach spaces. Isoperimetry and Processes.
Springer-Verlag.


\bibitem{SKgen} Panchenko, D. (2005). Free energy in the 
generalized Sherrington-Kirkpatrick mean field model.
{\it Rev. Math. Phys.} \textbf{17}, no. 7, 793-857.

\bibitem{PT} Parisi, G., Talagrand, M. (2004) 
On the distribution of the overlaps at given disorder.  
{\it C. R. Math. Acad. Sci. Paris,}  \textbf{339},  no. 4, 303-306.

\bibitem{SherK} Sherrington, D., Kirkpatrick, S. (1972)
Solvable model of a spin glass. 
{\it Phys. Rev. Lett.} {\bf 35}, 1792-1796.

\bibitem{SG} Talagrand, M. (2003)   
Spin Glasses: a Challenge for Mathematicians. 
Springer-Verlag.

\bibitem{T2} Talagrand, M. (2003) 
 On Guerra's broken replica-symmetry bound. 
{\it C. R. Math. Acad. Sci. Paris} {\bf 337}, no. 7, 477-480.

\bibitem{T-P} Talagrand, M. (2006)
Parisi formula. {\it Ann. of Math. (2)}  \textbf{163},  no. 1, 221-263.

\bibitem{TCortona} Talagrand, M. (2007)  
Large deviations, Guerra's and A.S.S. Schemes, and the Parisi hypothesis.
{\it Lecture Notes in Mathematics}, Vol. 1900, Eds: E. Bolthausen, A. Bovier.


\end{thebibliography}
\end{document}